\documentclass[12pt]{amsart}

 \usepackage{amsfonts,graphics,amsmath,amsthm,amsfonts,amscd,amssymb,amsmath,latexsym,multicol}
\usepackage{epsfig,url}
\usepackage{flafter}

\makeatletter

\def\jobis#1{FF\fi
  \def\predicate{#1}%
  \edef\predicate{\expandafter\strip@prefix\meaning\predicate}%
  \edef\job{\jobname}%
  \ifx\job\predicate
}

\makeatother

\if\jobis{proposal}%
\else
\fi

 \usepackage[matrix, arrow]{xy}

\DeclareMathOperator{\Supp}{Supp}

 \newcommand{\N}{\mathbb N}
 
 \newcommand{\Q}{\mathbb Q}
 \newcommand{\R}{\mathbb R}
 \newcommand{\Z}{\mathbb Z}
 \newcommand{\bir}{\dashrightarrow}
 \newcommand{\rddown}[1]{\left\lfloor{#1}\right\rfloor} 

 \numberwithin{equation}{subsection}
 \numberwithin{footnote}{subsection}

 \newtheorem{cor}[subsection]{Corollary}
 \newtheorem{lem}[subsection]{Lemma}
 \newtheorem{prop}[subsection]{Proposition}
 \newtheorem{thm}[subsection]{Theorem}
 \newtheorem{conj}[subsection]{Conjecture}

\newtheorem{claim}[subsection]{Claim}

{
    \newtheoremstyle{upright}%
        {8pt plus2pt minus4pt}%
        {8pt plus2pt minus4pt}%
        {\upshape}%
        {}%
        {\bfseries\scshape}%
        {}%
        {1em}%
        {}%

\theoremstyle{upright}

 \newtheorem{defn}[subsection]{Definition}
 
 \newtheorem{constr}[subsection]{Construction}

 \newtheorem{rem}[subsection]{Remark}

}

 \newcommand{\ke}[1]{$\acute{\mbox{e}}$}
 \newcommand{\ku}[1]{$\acute{\mbox{u}$}}
 \newcommand{\kl}[1]{$\acute{\mbox{l}}$}
 \newcommand{\kh}[1]{$\acute{\mbox{h}}$}
 \newcommand{\kr}[1]{$\acute{\mbox{r}}$}
 \newcommand{\kx}[1]{$\acute{\mbox{x}}$}
 \newcommand{\ki}[1]{${\^\i}$}

\title{On existence of log minimal models}
\author{Caucher Birkar}\thanks{2000 Mathematics Subject Classification: 14E30}
\date{\today}
\begin{document}
\maketitle

\begin{abstract}
In this paper, we prove that the log minimal model program in dimension $d-1$ implies the existence
of log minimal models for effective lc pairs (eg of nonnegative Kodaira dimension)
in dimension $d$. In fact, we prove that the same conclusion follows from a weaker assumption, namely, 
the log minimal model program  
with scaling in dimension $d-1$. This enables us to prove that effective lc pairs in dimension five have log minimal models. 
We also give new proofs of the existence of log minimal models for effective lc pairs in dimension four and 
of Shokurov reduction theorem.
\end{abstract}


\section{Introduction}

All the varieties in this paper are assumed to be over an algebraically closed
field $k$ of characteristic zero. See section $2$ for notations and terminology. For basic notions of the log minimal model 
program (LMMP) which are not specified below, eg singularities of pairs, we follow Koll\'ar-Mori [\ref{KM}] though we would 
 need their analogues for $\R$-divisors.

One of the main problems in birational geometry and the classification theory of algebraic varieties in the last three decades 
or so has been the following 

\begin{conj}[Minimal model]\label{mmodel}
Any lc pair $(X/Z,B)$ has a log minimal model or a Mori fibre space.
\end{conj}

Here $B$ is an $\R$-boundary. Roughly speaking, a model $Y$ birational to $X$ is a log minimal model if $K+B$ is  
nef on it, and a Mori fibre space if it has a log Fano fibre structure 
which is negative with respect to $K+B$ where $K$ stands for the canonical divisor. 

The two-dimensional case of the above conjecture is considered to be a classical early 20th 
century result of the Italian algebraic geometry school, 
at least when $X$ is smooth and $B=0$. The three-dimensional case was proved by the 
contributions of many people, in particular, Mori's theorems on extremal rays [\ref{Mori-notnef}] and existence of flips 
[\ref{Mori-flip}], Shokurov's results on existence of log flips [\ref{log-flips}], termination [\ref{nonvanishin}][\ref{log-models}],  
and nonvanishing [\ref{nonvanishin}], 
Kawamata-Viehweg vanishing theorem [\ref{Kaw-vanishing}][\ref{Viehweg-vanishing}], 
Shokurov-Kawamata base point free theorem [\ref{nonvanishin}], and Kawamata's termination of log flips [\ref{Kaw-termination}]. 
A much simpler proof of Conjecture \ref{mmodel} in dimension three would be a combination of Shokurov's simple proof of 
existence of log flips in dimension three [\ref{pl-flips}], his termination in the terminal case [\ref{nonvanishin}] and  his recent method of constructing 
log minimal models [\ref{ordered}](see also [\ref{B2}]). In case the pair is effective, i.e. there is an effective $\R$-divisor $M\equiv K_X+B/Z$, 
Shokurov's existence of log flips in dimension three [\ref{pl-flips}] 
and this paper give yet another proof.

In dimension four, Shokurov's existence of log flips [\ref{pl-flips}] and Kawamata-Matsuda-Matsuki termination theorem [\ref{KMM-mmp}] 
prove the conjecture when $B=0$ and $X$ has terminal singularities. The general case in dimension four is a theorem of 
Shokurov [\ref{ordered}] with a short proof given in Birkar [\ref{B2}] in the klt case. In dimension five, the conjecture 
is proved below when $(X/Z,B)$ is effective, in particular, 
if the pair has nonnegative Kodaira dimension. 
In higher dimension, the conjecture is known in the klt case when $B$ is big$/Z$ by 
Birkar-Cascini-Hacon-M$^c$Kernan [\ref{BCHM}], in particular, if $X/Z$ is of general type or a flipping contraction.

To construct a log minimal model or a Mori fibre space one often runs the LMMP. In order to be able to run the 
LMMP we have to deal with the existence and termination of log flips. 
As experience teaches us, the best way might be an inductive argument. Existence of log flips 
in the klt (and $\Q$-factorial dlt) case is treated in [\ref{BCHM}] and [\ref{BP}] but we will not use their results below since we want to 
consider existence of log flips as a special case of existence of log minimal models hence to fit it 
into the inductive approach below and the results of Shokurov [\ref{pl-flips}] and Hacon-M$^c$Kernan [\ref{HM2}]. 
As for the termination, inductive arguments 
involve other conjectures:  the LMMP in dimension $d-1$ and the ascending chain condition 
(ACC) for lc thresholds in dimension $d$ imply termination of log flips for effective lc pairs in dimension 
$d$ by Birkar [\ref{B}]. Other results in this direction include 
the reduction of termination of log flips in dimension $d$ to  the ACC and semi-continuity for minimal log discrepancies 
in dimension $d$ by Shokurov [\ref{mld's}], and the reduction of termination of log flips for effective lc pairs in dimension $d$ to 
LMMP, boundedness of certain Fano varieties and 
the ACC for minimal log discrepancies in dimension $d-1$ by Birkar-Shokurov [\ref{BSh}]. 
 
In this paper, by focusing on constructing log minimal models rather than proving a general termination we 
prove the following inductive

\begin{thm}\label{induction}
Assume the LMMP for $\Q$-factorial dlt pairs of dimension $d-1$. Then, Conjecture \ref{mmodel} holds 
for effective lc pairs in dimension $d$.
\end{thm}

We must emphasize that unlike in [\ref{BCHM}] and [\ref{BP}], we do not assume bigness of the boundary $B$. 
The abundance conjecture predicts that a lc pair $(X/Z,B)$ has a log minimal model exactly when it is 
effective. When $(X/Z,B)$  is not effective 
existence of a Mori fibre space is predicted. 

In Theorem \ref{induction}, the full LMMP in dimension $d-1$ is not necessary; we only need LMMP with scaling (see Definition \ref{d-scaling}).

\begin{thm}\label{c-induction}
 Assume the LMMP with scaling for $\Q$-factorial dlt pairs in dimension $d-1$. Then, Conjecture \ref{mmodel} 
holds for effective lc pairs in dimension $d$.
\end{thm}

In a forthcoming paper, we reduce Conjecture \ref{mmodel} to the following
weak nonvanishing conjecture. 

\begin{conj}\label{effectivity}
If a $\Q$-factorial dlt pair $(X/Z,B)$ is pseudo-effective, then it is effective.
\end{conj}

  When $d=5$,
instead of proving Conjecture \ref{effectivity} in dimension four we directly use termination of
log flips with scaling in dimension four which we prove that it follows from
results of Shokurov [\ref{ordered}] and Alexeev-Hacon-Kawamata [\ref{AHK}].

\begin{cor}\label{5-fold}
Log minimal models exist for effective lc pairs in dimension five.
\end{cor}

We also get a new proof of the following result which was first proved by Shokurov [\ref{ordered}] (see also [\ref{B2}]).

\begin{cor}\label{4-fold}
Log minimal models exist for effective lc pairs in dimension four.
\end{cor}

Our method gives a new proof of Shokurov
reduction theorem, that is, reducing the existence of log flips to the special termination and the existence of pl flips 
(see Theorem \ref{t-reduction}). 
A variant of the reduction theorem is  an important ingredient in the construction of log flips in [\ref{BCHM}][\ref{BP}]. 

For other important applications see section two of [\ref{BP}] and section five of [\ref{BCHM}]. 
To avoid confusion, let us make it clear that this paper is logically independent of [\ref{BP}] and [\ref{BCHM}]

\section*{Acknowledgement}

 I would like to thank Yoichi Miyaoka and Tokyo University, where part of this work was done in March 2007, for their support and
 hospitality. I also thank Stephane Druel, Christopher Hacon and James M$^c$Kernan for their comments.


\section{Basics}

Let $k$ be an algebraically closed field of characteristic zero fixed throughout the paper. 

A \emph{pair} $(X/Z,B)$ consists of normal quasi-projective varieties $X,Z$ over $k$, an $\R$- divisor $B$ on $X$ with
coefficients in $[0,1]$ such that $K_X+B$ is $\mathbb{R}$-Cartier, and a projective 
morphism $X\to Z$. For a prime divisor $D$ on some birational model of $X$ with a
nonempty centre on $X$, $a(D,X,B)$
denotes the log discrepancy.

A pair $(X/Z,B)$ is called \emph{pseudo-effective} if $K_X+B$ is pseudo-effective/$Z$, that is,
up to numerical equivalence/$Z$ it is the limit of effective $\R$-divisors. The pair is called \emph{effective} if
$K_X+B$ is effective/$Z$, that is, there is an $\R$-divisor $M\ge 0$ such that $K_X+B\equiv M/Z$;
in this case, we call $(X/Z,B,M)$ a \emph{triple}. By a log resolution of a triple $(X/Z,B,M)$
we mean a log resolution of $(X, \Supp B+M)$. When we refer to a triple as being lc, dlt, etc, we mean that 
the underlying pair has such properties.

Let $(X/Z,B)$ be a lc pair. By a \emph{log flip}$/Z$ we mean the flip of a $K_X+B$-negative extremal flipping contraction$/Z$ [\ref{B}, Definition 2.3], 
and by a \emph{pl flip}$/Z$ we mean a log flip$/Z$ such that $(X/Z,B)$ is $\Q$-factorial dlt and the log flip is also an $S$-flip for 
some component $S$ of $\rddown{B}$.

A \emph{sequence of log flips$/Z$ starting with} $(X/Z,B)$ is a sequence $X_i\bir X_{i+1}/Z_i$ in which  
$X_i\to Z_i \leftarrow X_{i+1}$ is a $K_{X_i}+B_i$-flip$/Z$, $B_i$ is the birational transform 
of $B_1$ on $X_1$, and $(X_1/Z,B_1)=(X/Z,B)$.

In this paper, \emph{special termination} means termination near $\rddown{B}$ of any sequence of log flips$/Z$ 
starting with a pair $(X/Z,B)$, that is, 
the log flips do not intersect $\rddown{B}$ after finitely many of them. There is a more general notion of special termination 
which claims that the log 
flips do not intersect any lc centre after finitely many steps but we do not use it below. 

\begin{defn}
For an $\R$-divisor $D=\sum d_iD_i$, let $D^{\le 1}:=\sum d_i'D_i$ where $d_i'=\min\{d_i,1\}$. As usual, $D_i$ are distinct 
prime divisors.
\end{defn}

\begin{defn}
For a triple  $(X/Z,B,M)$,
define
$$
\theta(X/Z,B,M):=\#\{i~|~m_i\neq 0 ~~\mbox{and}~~ b_i\neq 1\}
$$
where $B=\sum b_iD_i$ and $M=\sum m_iD_i$.
\end{defn}

\begin{constr}\label{smooth-model}
Let $(X/Z,B,M)$ be a lc triple. Let $f\colon W\to X$ be a log resolution of
$(X/Z,B,M)$, and let
$$
B_W:=B^{\sim}+\sum E_j
$$
where $\sim$ stands for the birational transform, and $E_j$ are the prime exceptional divisors of $f$.

Obviously, $(W/Z,B_W)$ is $\Q$-factorial dlt, and it is effective because
$$
M_W :=K_W+B_W-{f^*(K_X+B)}+{f^*M} \equiv K_W+B_W/Z
$$
is effective. Note that since $(X/Z,B)$ is lc,
$$
K_W+B_W-{f^*(K_X+B)}\ge 0
$$
In addition, each component of
$M_W$ is either a component of $M^{\sim}$ or an exceptional divisor $E_j$. Thus, by construction
$$
\theta(W/Z,B_W,M_W)=\theta(X/Z,B,M)
$$

We call $(W/Z,B_W,M_W)$ and $(W/Z,B_W)$  \emph{log smooth 
models} of $(X/Z,B,M)$ and $(X/Z,B)$ respectively.
\end{constr}

\begin{defn}[Cf. {[\ref{log-flips}]}]\label{d-model}
A pair $(Y/Z,B+E)$ is a \emph{log birational model} of $(X/Z,B)$ if we have a birational map
$\phi\colon X\bir Y/Z$, $B$ on $Y$ which is the birational transform of $B$ on $X$
(for simplicity we use the same notation), and $E=\sum E_j$ where $E_j$ are the
exceptional/$X$ prime divisors of $Y$. $(Y/Z,B+E)$ is a \emph{nef model} of $(X/Z,B)$ if in addition\\\\
(1) $(Y/Z,B+E)$ is $\Q$-factorial dlt, and\\
(2) $K_Y+B+E$ is nef/$Z$.\\

And  we call $(Y/Z,B+E)$ a \emph{log minimal model} of $(X/Z,B)$ if in addition\\\\
(3) for any prime divisor $D$ on $X$ which is exceptional/$Y$, we have
$$
a(D,X,B)<a(D,Y,B+E)
$$
\end{defn}

\begin{defn}[Mori fibre space]
A log birational model $(Y/Z,B+E)$ of a lc pair $(X/Z,B)$ is called a Mori fibre space if 
$(Y/Z,B+E)$ is $\Q$-factorial dlt, there is a $K_Y+B+E$-negative extremal contraction $Y\to T/Z$ 
with $\dim Y>\dim T$, and 
$$
a(D,X,B)\le a(D,Y,B+E)
$$
for  any prime divisor $D$ (on birational models of $X$) and the strict inequality holds if $D$ is on $X$ and contracted$/Y$.
\end{defn}

Our definitions of log minimal models and Mori fibre spaces are slightly different 
from those in [\ref{KM}], the difference being that we do not assume that $\phi^{-1}$ does not contract divisors. 
However, in the plt case, our definition of log minimal models and that of [\ref{KM}] coincide (see Remark \ref{r-mmodels} (iii)).

\begin{rem}\label{r-mmodels}
Let $(X/Z,B)$ be a lc pair.\\

(i) Suppose that  $(W/Z,B_W)$ is a log smooth model of $(X/Z,B)$ and $(Y/Z,B_W+E)$ 
a log minimal model of 
$(W/Z,B_W)$. We can also write $(Y/Z,B_W+E)$ as $(Y/Z,B+E')$ where $B$ on $Y$ is the birational transform of 
$B$ on $X$ and $E'$ is the reduced divisor whose components are the exceptional$/X$ divisors on $Y$. 
 Let $D$ be a prime divisor on $X$ contracted$/Y$. Then, 
$$
a(D,X,B)=a(D,W,B_W)<a(D,Y,B_W+E)=a(D,Y,B+E')
$$
which implies that $(Y/Z,B+E')$ is a log minimal model of $(X/Z,B)$.\\
 
(ii) Let $(Y/Z,B+E)$ be a log minimal model of $(X/Z,B)$ 
and take a common resolution $f\colon W\to X$ 
and $g\colon W\to Y$, then 
$$
f^*(K_X+B)\ge g^*(K_Y+B+E)
$$
 by applying the negativity lemma [\ref{log-flips}, 1.1].
This, in particular, means that $a(D,X,B)\le a(D,Y,B+E)$ for any prime divisor $D$ (on birational models of $X$). 
Moreover, if $(X/Z,B,M)$ is a triple, then $(Y/Z,B+E, g_*f^*M)$ is also a triple.\\

(iii) If $(X/Z,B)$ is plt and $(Y/Z,B+E)$ is a log minimal model, for any component $D$ of $E$ on $Y$, 
$0<a(D,X,B)\le a(D,Y,B+E)=0$ which is not possible, so $E=0$.
\end{rem}

\section{Proofs}

Before getting into the proof of our results we need some preparation.

\begin{lem}\label{ray}
Let $(X/Z,B+C)$ be a $\Q$-factorial lc pair where $B,C\ge 0$, 
$K_X+B+C$ is nef/$Z$, and $(X/Z,B)$ is dlt. Then, either $K_X+B$ is also nef/$Z$ or there is an extremal ray $R/Z$ such
that $(K_X+B)\cdot R<0$, $(K_X+B+\lambda C)\cdot R=0$, and $K_X+B+\lambda C$ is nef$/Z$ where
$$
\lambda:=\inf \{t\ge 0~|~K_X+B+tC~~\mbox{is nef/$Z$}\}
$$
\end{lem}
\begin{proof}
Suppose that $K_X+B$ is not nef$/Z$ and let $\{R_i\}_{i\in I}$ be the set of $(K_X+B)$-negative extremal rays/$Z$
and $\Gamma_i$ an extremal curve of $R_i$ [\ref{ordered}, Definition 1]. Let $\mu:=\sup \{\mu_i\}$ where
$$
\mu_i:=\frac{-(K_X+B)\cdot \Gamma_i}{C\cdot \Gamma_i}
$$

Obviously, $\lambda=\mu$ and $\mu\in (0,1]$. It is enough to
prove that $\mu=\mu_l$ for some $l$. 
By [\ref{ordered}, Proposition 1], there are positive real numbers
$r_1,\dots, r_s$ and a positive integer $m$ (all independent of $i$) such that
$$
(K_X+B)\cdot \Gamma_i=\sum_{j=1}^s\frac{r_jn_{i,j}}{m}
$$
where $-2(\dim X)m\le n_{i,j}\in\Z$.
On the other hand, by [\ref{log-models}, First Main Theorem 6.2, Remark 6.4] we can write 
$$
K_X+B+C=\sum_{k=1}^{t} r_k'(K_X+\Delta_k)
$$ 
where $r_1',\cdots, r_{t}'$ are positive 
real numbers such that for any $k$ we have: $(X/Z,\Delta_k)$ is lc with $\Delta_k$ being rational, and $(K_X+\Delta_k)\cdot \Gamma_i\ge 0$ for any $i$. Therefore, 
there is a positive integer $m'$ (independent of $i$) such that 
$$
(K_X+B+C)\cdot \Gamma_i=\sum_{k=1}^{t}\frac{r_k'n_{i,k}'}{m'}
$$
where $0\le n_{i,k}'\in\Z$. 

The set $\{n_{i,j}\}_{i,j}$ is finite.
Moreover,
$$
\frac{1}{\mu_i}= \frac{C\cdot \Gamma_i}{-(K_X+B)\cdot \Gamma_i}=\frac{(K_X+B+C)\cdot \Gamma_i}{-(K_X+B)\cdot \Gamma_i}+1
=-\frac{m\sum_k r_k'n_{i,k}'}{m'\sum_j r_jn_{i,j}}+1
$$

Thus, $\inf \{\frac{1}{\mu_i}\}=\frac{1}{\mu_l}$ for some $l$ and so $\mu=\mu_l$.
\end{proof}

\begin{defn}[LMMP with scaling]\label{d-scaling}
Let $(X/Z,B+C)$ be a lc pair such that $K_X+B+C$ is nef/$Z$, $B\ge 0$, and $C\ge 0$ is $\R$-Cartier. 
Suppose that either $K_X+B$ is nef/$Z$ or there is an extremal ray $R/Z$ such
that $(K_X+B)\cdot R<0$, $(K_X+B+\lambda_1 C)\cdot R=0$, and $K_X+B+\lambda_1 C$ is nef$/Z$ where
$$
\lambda_1:=\inf \{t\ge 0~|~K_X+B+tC~~\mbox{is nef/$Z$}\}
$$
 When $(X/Z,B)$ is $\Q$-factorial dlt, the last sentence follows from 
Lemma \ref{ray}. If $R$ defines a Mori fibre structure, we stop. Otherwise assume that $R$ gives a divisorial 
contraction or a log flip $X\bir X'$. We can now consider $(X'/Z,B'+\lambda_1 C')$  where $B'+\lambda_1 C'$ is 
the birational transform 
of $B+\lambda_1 C$ and continue the argument. That is, suppose that either $K_{X'}+B'$ is nef/$Z$ or 
there is an extremal ray $R'/Z$ such
that $(K_{X'}+B')\cdot R'<0$, $(K_{X'}+B'+\lambda_2 C')\cdot R'=0$, and $K_{X'}+B'+\lambda_2 C'$ is nef$/Z$ where
$$
\lambda_2:=\inf \{t\ge 0~|~K_{X'}+B'+tC'~~\mbox{is nef/$Z$}\}
$$
 By continuing this process, we obtain a 
special kind of LMMP$/Z$ which is called the \emph{LMMP$/Z$ on $K_X+B$ with scaling of $C$}; note that it is not unique. 
This kind of LMMP was first used by Shokurov [\ref{log-flips}].
When we refer to \emph{termination with scaling} we mean termination of such an LMMP.

\emph{Special termination  with scaling} means termination near $\rddown{B}$ of any sequence of log flips$/Z$ with scaling 
of $C$, i.e. after finitely many steps, the locus of the extremal rays in the process do not intersect $\Supp \rddown{B}$.

When we have a lc pair $(X/Z,B)$, we can always find an ample$/Z$ $\R$-Cartier divisor $C\ge 0$ such that 
$K_X+B+C$ is lc and nef$/Z$,  so we can run the LMMP$/Z$ with scaling assuming that all the 
necessary ingredients exist, eg extremal rays, log flips. 
\end{defn}

\begin{lem}\label{l-extraction}
 Assume the special termination with scaling for $\Q$-factorial dlt pairs in dimension $d$ 
and the existence of pl flips in dimension $d$. Let $(X/Z,B)$ be a lc pair of dimension
 $d$ and let $\{D_i\}_{i\in I}$ be a finite set of exceptional$/X$ prime divisors (on birational
 models of $X$) such that $a(D_i,X,B)\le 1$. Then, there is a $\Q$-factorial dlt pair $(Y/X,B_Y)$
 such that\\\\
 (1) $Y/X$ is birational and $K_Y+B_Y$ is the crepant pullback of $K_X+B$,\\
 (2) every exceptional/$X$ prime divisor $E$ of $Y$ is one of the $D_i$ or $a(E,X,B)=0$,\\
 (3) the set of exceptional/$X$ prime divisors of $Y$ includes $\{D_i\}_{i\in I}$.\\
\end{lem}
\begin{proof}
 Let $f\colon W\to X$ be a log resolution of
$(X/Z,B)$ and let $\{E_j\}_{j\in J}$
be the set of prime exceptional divisors of $f$. Moreover, we can assume that
for some $J'\subseteq J$, $\{E_j\}_{j\in J'}=\{D_i\}_{i\in I}$.
Now define
$$
K_W+B_W:=f^*(K_{X}+B)+\sum_{j\notin J'} a(E_j,X,B)E_j
$$
which ensures that if $j\notin J'$, then $E_j$ is a component of $\rddown{B}$.
 
By running the LMMP/$X$ on $K_W+B_W$ with scaling of a suitable ample$/X$ $\R$-divisor, and 
using the special termination we get a log minimal model of $(W/X,B_W)$ which we may denote by  
$(Y/X,B_Y)$. By the negativity lemma, all the $E_j$ are contracted in the
process except if $j\in J'$ or if $a(E_j,X,B)=0$. By construction, $K_Y+B_Y$ is the crepant pullback of  $K_X+B$.
\end{proof}


\begin{prop}\label{p-main}
Assume the special termination with scaling for $\Q$-factorial dlt pairs in dimension $d$ and 
the existence of pl flips in dimension $d$. Then, any effective lc pair in dimension $d$ has a 
log minimal model.
\end{prop}

\begin{proof}
\emph{Step 1.} Let $\mathfrak W$ be the set of triples $(X/Z,B,M)$ such that\\

(1) $(X/Z,B)$ is lc of dimension $d$,

(2) $(X/Z,B)$ does not have a log minimal model.\\

It is enough to prove that $\mathfrak W$ is empty. Assume otherwise and
choose  $(X/Z,B,M)\in\mathfrak W$ with minimal $\theta(X/Z,B,M)$.
Replace $(X/Z,B,M)$ with a log smooth model as in Construction
\ref{smooth-model} which preserves the minimality of $\theta(X/Z,B,M)$.

If $\theta(X/Z,B,M)=0$, then either $M=0$ in which case we already have a log minimal model,
or by running the LMMP/$Z$ on $K_X+B$ with scaling of a suitable ample$/Z$ $\R$-divisor we get a log 
minimal model because by the special termination,
flips and divisorial contractions will not intersect $\Supp \rddown{B}\supseteq \Supp M$ after finitely many steps.
This is a contradiction. 
Note that we only need pl flips here. We may then assume that $\theta(X/Z,B,M)> 0$.\\

\emph{Step 2.} Define
$$
\alpha:=\min\{t>0~|~~\rddown{(B+tM)^{\le 1}}\neq \rddown{B}~\}
$$

In particular, $(B+\alpha M)^{\le 1}=B+C$ for some $C\ge 0$ supported on $\Supp M$, and $\alpha M=C+M'$ 
where $M'$ is supported on $\Supp \rddown{B}$. Thus, outside $\Supp \rddown{B}$ 
we have $C=\alpha M$. The pair $(X/Z,B+C)$ is  $\Q$-factorial dlt and
$(X/Z,B+C,M+C)$ is a triple. By construction
$$
\theta(X/Z,B+C,M+C)<\theta(X/Z,B,M)
$$
and so $(X/Z,B+C,M+C)\notin \mathfrak W$. Therefore, $(X/Z,B+C)$ has a log minimal model,
say $(Y/Z,B+C+E)$. By definition,
$
K_Y+B+C+E
$
is nef/$Z$.\\

\emph{Step 3.} Now run the LMMP$/Z$ on $K_Y+B+E$ with scaling of $C$. Note that we only need pl flips here 
because every extremal ray contracted would have negative intersection with some component of $\rddown{B}+E$ 
by Remark \ref{r-mmodels} (ii) and the properties of $C$ mentioned in Step 2. 
By the special termination,
after finitely many steps, $\rddown{B}+E$ does not intersect the extremal rays
contracted by the LMMP hence we end up with a model $Y'$ on which $K_{Y'}+B+E$ is nef/$Z$. Clearly, $(Y'/Z,B+E)$
is a nef model of  $(X/Z,B)$ but may not be a log minimal model because condition (3) of Definition
\ref{d-model} may not be satisfied.\\

\emph{Step 4.} Let
$$
\mathcal T=\{t\in [0,1]~|~ (X/Z,B+tC)~~\mbox{has a log minimal model}\}
$$
Since $1\in\mathcal T$, $\mathcal T\neq \emptyset$. Let $t\in\mathcal T\cap (0,1]$ and let 
$(Y_t/Z,B+tC+E)$ be any log minimal model of $(X/Z,B+tC)$. Running the LMMP/$Z$ on $K_{Y_t}+B+E$ with scaling
of $tC$ shows that there is $t'\in(0,t)$ such that $[t',t]\subset \mathcal{T}$
because condition (3) of Definition \ref{d-model} is an open condition. The
LMMP terminates for the same reasons as in Step 3 and we note again that the log flips 
required are all pl flips.\\

\emph{Step 5.} Let $\tau=\inf \mathcal T$. If $\tau\in\mathcal{T}$, then by Step 4, $\tau=0$
and so we are done by deriving a contradiction. Thus, we may assume that $\tau\notin\mathcal{T}$. In this case, there is a sequence
$t_1>t_2>\cdots$ in $ \mathcal{T}\cap (\tau,1]$
such that $\lim t_k=\tau$. For each $t_k$ let $(Y_{t_k}/Z,B+t_kC+E)$ be any log minimal model of
$(X/Z,B+t_kC)$ which exists by definition of $\mathcal{T}$ and from which we get a nef model $(Y_{t_k}'/Z,B+\tau C+E)$
for $(X/Z,B+\tau C)$ by running the LMMP/$Z$ on $K_{Y_{t_k}}+B+E$ with scaling of $t_kC$.
Let $D\subset X$ be a prime divisor contracted/$Y_{t_k}'$. If $D$ is contracted/$Y_{t_k}$, then 
$$
a(D,X,B+t_k C) < a(D,Y_{t_k},B+t_k C+E)
$$
$$
\le a(D,Y_{t_k},B+\tau C+E)\le a(D,Y_{t_k}',B+\tau C+E)
$$
 but if $D$ is not contracted/$Y_{t_k}$ we have
$$
a(D,X,B+t_k C) =a(D,Y_{t_k},B+t_k C+E)
$$
$$
\le a(D,Y_{t_k},B+\tau C+E) <a(D,Y_{t_k}',B+\tau C+E)
$$
because $(Y_{t_k}/Z,B+t_kC+E)$ is a log minimal model of $(X/Z,B+t_kC)$ and  $(Y_{t_k}'/Z,B+\tau C+E)$ 
is a log minimal model of $(Y_{t_k}/Z,B+\tau C+E)$. Thus, in any case we have 
$$
a(D,X,B+t_k C) < a(D,Y_{t_k}',B+\tau C+E)
$$

Replacing the sequence $\{t_k\}_{k\in \N}$ with a subsequence, we can assume that all the induced rational maps $X\bir Y_{t_k}'$
contract the same components of $B+\tau C$. 

\begin{claim}
$(Y_{t_k}'/Z,B+\tau C+E)$ and $(Y_{t_{k+1}}'/Z,B+\tau C+E)$ have equal log discrepancies for any $k$.
\end{claim}
\begin{proof}
Let $f\colon W\to Y_{t_k}'$ and  $g\colon W\to Y_{t_{k+1}}'$ be resolutions with a common $W$.
Put 
$$
F=f^*(K_{Y_{t_{k}}'}+B+\tau C+E)-g^*(K_{Y_{t_{k+1}}'}+B+\tau C+E)
$$
which is obviously anti-nef$/Y_{t_k}'$. So, by the negativity lemma, $F\ge 0$ iff $f_*F\ge 0$. 
 Suppose that $D$ is a component of $f_*F$ with negative coefficient. By assumptions, 
$D$ is not a component of $B+\tau C+E$ on $Y_{t_k}'$ and it must be exceptional$/Y_{t_{k+1}}'$. 
Moreover, the coefficient of $D$ in $f_*F$ is equal to 
$a(D,Y_{t_{k+1}}',B+\tau C+E)-1$. On the other hand, 
$$
1=a(D,X,B+t_{k+1} C) <a(D,Y_{t_{k+1}}',B+\tau C+E)
$$
a contradiction. On the other hand, $F$ is nef$/Y_{t_{k+1}}'$ and a similar argument would imply that $-F\ge 0$ 
hence $F=0$ and the claim 
follows. 
\end{proof}

Therefore, each $(Y_{t_k}'/Z,B+\tau C+E)$ is a nef model of $(X/Z,B+\tau C)$ such that
$$
 a(D,X,B+\tau C)=\lim a(D,X,B+t_k C)\leq a(D,Y_{t_k}',B+\tau C+E)
$$
for any prime divisor $D\subset X$ contracted/$Y_{t_k}'$.\\

\emph{Step 6.} To get a log minimal model of  $(X/Z,B+\tau C)$ we just need to extract those prime divisors $D$ on $X$ 
contracted$/Y_{t_{k}}'$ for which 
$$
 a(D,X,B+\tau C)=a(D,Y_{t_k}',B+\tau C+E)
$$

This is achieved by applying Lemma \ref{l-extraction} to construct a suitable crepant model
of $(Y_{t_k}',B+\tau C+E)$ which would then be a log minimal model of $(X/Z,B+\tau C)$.
Thus, $\tau\in\mathcal T$ and this gives a contradiction. Therefore, $\mathfrak{W}=\emptyset$.
\end{proof}

\begin{lem}\label{l-st}
LMMP with scaling in dimension $d-1$ for $\Q$-factorial dlt pairs implies the 
special termination with scaling in dimension $d$ for dlt pairs.
\end{lem}
\begin{proof}
Let $(X/Z,B+C)$ be a dlt pair such that $K_X+B+C$ is nef/$Z$, $B\ge 0$, and $C\ge 0$ is $\R$-Cartier. Suppose that 
we get a sequence $X_i\bir X_{i+1}/Z_i$ of log flips while running the LMMP$/Z$ with scaling 
of $C$. We may assume that $X=X_1$. Let $S$ be 
a component of $\rddown{B}$ and let $S_i$ be its birational transform on $X_i$ which is normal, and $T_i$ the normalisation 
of its birational transform on $Z_i$. Also put $K_{S_i}+B_{S_i}:=(K_{X_i}+B_i)|_{S_i}$ and $C_{S_i}=C_i|_{S_i}$ 
where $B_i$ and $C_i$ on $X_i$ are the birational transforms of $B$ and $C$.

 In general,  the induced birational map $S_i\bir S_{i+1}$ is not a log flip, however, if $i\gg 0$ 
it is an isomorphism in codimension one; it could also be an isomorphism. Moreover, 
$S_i\to T_i$ is not necessarily extremal. Now pick $i\gg 0$. Since $K_{S_i}+B_{S_i}$ is dlt, 
by the proof of Lemma \ref{l-extraction}, we can construct a $\Q$-factorialisation of 
$S_i$, say $S_i'$ such that $S_i'\to S_i$ is small. Assume that $K_{S_i'}+B_{S_i'}$ and $C_{S_i'}$ are the 
pullbacks of $K_{S_i}+B_{S_i}$ and $C_{S_i}$. By assumptions, after a finite 
 sequence of $K_{S_i'}+B_{S_i'}$-log flips$/T_i$ we end up with $S_{i+1}'/S_{i+1}$ on which 
$K_{S_{i+1}'}+B_{S_{i+1}'}$ the pullback of $K_{S_{i+1}}+B_{S_{i+1}}$ is nef$/T_i$. Since  
$K_{S_i'}+B_{S_i'}+\lambda_iC_{S_i'}\equiv 0/T_i$ the sequence can also be considered as a sequence of 
log flips with scaling of $\lambda_iC_{S_i'}$. 
Continuing the argument gives a sequence of log flips$/Z$ with scaling which must terminate by assumptions. 
So, the sequence  $X_i\bir X_{i+1}/Z_i$ does not intersect $S$ after finitely many log flips.
\end{proof}

\begin{rem}\label{r-pl-flips}
Assume the LMMP with scaling for $\Q$-factorial dlt pairs in dimension 
$d-1$, then pl flips exist in dimension $d$ by [\ref{HM2}].
\end{rem}

\begin{proof}(of Theorem \ref{c-induction}).
Immediate by Lemma \ref{l-st}, Remark \ref{r-pl-flips}, and Proposition \ref{p-main}.
\end{proof}

\begin{proof}(of Theorem \ref{induction}).
Immediate by Theorem \ref{c-induction}.
\end{proof}

We now turn to the proof of the corollaries.

\begin{lem}\label{st}
LMMP with scaling holds for $\Q$-factorial dlt pairs in dimension $4$.
\end{lem}
\begin{proof}
Log flips exist by [\ref{pl-flips}][\ref{HM2}] so it is enough to verify the termination. 
Let $(X/Z,B+C)$ be a $\Q$-factorial dlt pair of dimension $4$. 
Suppose that we get a sequence of log flips $X_i\bir X_{i+1}/Z_i$, when running the LMMP$/Z$ with scaling
of $C$, which does not terminate. We may assume that $X=X_1$. Let $\lambda=\lim \lambda_i$  where the $\lambda_i$
are obtained as in Definition \ref{d-scaling} for this LMMP/$Z$.
By applying [\ref{ordered}, Corollary 12] we deduce that $\lambda_i$ 
stabilize, that is, $\lambda=\lambda_i$
for $i\gg 0$. Therefore, we get an infinite sequence of $K_X+B+\frac{1}{2}\lambda C$-flips
such that $C$ is positive on each of these flips. 
Now by special termination in dimension $4$ and by [\ref{AHK}, Theorem 2.15], 
we may assume that each of the log flips is of type $(2,1)$, that is, 
the flipping locus is of dimension $2$ and the flipped locus is of dimension $1$. 
Take the normalisation of a  compactification of $X=X_1$ and denote it by $Y_1$. 
This also induces compatible compactification 
and normalisation for each $X_i$ denoted by $Y_i$; we get a  sequence $Y_i\bir Y_{i+1}$ of birational maps which 
are isomorphisms in codimension one and such that 
a codimension two cycle disappears under each  $Y_i\bir Y_{i+1}$. Now we get a contradiction by [\ref{mld's}, Lemma 1] 
 essentially 
for topological reasons.
\end{proof}

\begin{proof}(of Corollary \ref{5-fold})
Immediate by Lemma \ref{st} and Theorem \ref{c-induction}.
\end{proof}

\begin{proof}(of Corollary \ref{4-fold})
Immediate by [\ref{log-models}] and Theorem \ref{c-induction}.
\end{proof}

Next we show that Shokurov reduction theorem [\ref{pl-flips}, Reduction Theorem 1.2] follows easily from Proposition \ref{p-main}. 

\begin{thm}[Shokurov reduction]\label{t-reduction}
Assume the special termination for $\Q$-factorial dlt pairs in dimension $d$ and the existence of pl flips in dimension $d$.
Then, log flips exist for klt (hence $\Q$-factorial dlt) pairs in dimension $d$. 
\end{thm}
\begin{proof}
Let $(X/Z,B)$ be a klt pair of dimension $d$ and $f\colon X\to Z'$  a
$(K_X+B)$-flipping contraction$/Z$. We can apply Proposition \ref{p-main} to construct a
log minimal model $(Y/Z',B+E)$ of $(X/Z',B)$.  Now since $(X/Z',B)$ is klt,
$E=0$. So, $(Y/Z',B)$ is also klt and by the base point free theorem [\ref{HM2}, Theorem 5.2.1] 
it has a log canonical model which gives the flip of $f$. 
\end{proof}

As mentioned in Remark \ref{r-pl-flips}, by [\ref{HM2}], the LMMP with scaling for $\Q$-factorial dlt 
pairs in dimension $d-1$ implies the existence of pl flips in dimension $d$. In Lemma \ref{l-st}, we proved 
that the LMMP with scaling for $\Q$-factorial dlt pairs in dimension $d-1$ implies special termination 
with scaling in dimension $d$ for dlt pairs. Thus, the previous 
theorem can be restated as saying that the LMMP with scaling for $\Q$-factorial dlt pairs in dimension $d-1$ 
implies the existence of log flips for klt (hence $\Q$-factorial dlt) pairs.

\vspace{2cm}

\flushleft{DPMMS}, Centre for Mathematical Sciences,\\
Cambridge University,\\
Wilberforce Road,\\
Cambridge, CB3 0WB,\\
UK\\
email: c.birkar@dpmms.cam.ac.uk

\end{document}